\documentclass[11pt]{amsart}
\usepackage{amsfonts}
\usepackage{amsmath}
\usepackage{amscd}
\usepackage{amsthm}
\usepackage{enumerate}

\usepackage{graphicx}

\usepackage{color}

\usepackage{hyperref}

\usepackage[margin=0.86in]{geometry}

\newtheorem{theorem}{Theorem} 
\newtheorem*{theorem*}{Theorem}

\newtheorem{proposition}[theorem]{Proposition}

\newtheorem{corollary}[theorem]{Corollary}

\theoremstyle{remark}
\newtheorem{rmk}[theorem]{Remark}
\newtheorem{asmp}[theorem]{Assumption}
\newtheorem{example}[theorem]{Example}

\newcommand{\pp}{\mathbb{P}}
\newcommand{\qq}{\mathbb{Q}}
\newcommand{\rr}{\mathbb{R}}

\newcommand{\nn}{\mathbb{N}}
\newcommand{\ep}{\hfill \ensuremath{\Box}}
\newcommand{\eq}{\begin{equation}}
\newcommand{\en}{\end{equation}}

\newcommand{\ev}{\mathbb{E}}

\numberwithin{equation}{section} \numberwithin{theorem}{section}

\title[Infinite particle systems]{A construction of infinite Brownian particle systems}

\author{Mykhaylo Shkolnikov}
\address{Department of Mathematics, Princeton University, Princeton, NJ 08544, USA}
\email{mshkolni@gmail.com}

\begin{document}

\begin{abstract}
The paper identifies families of quasi-stationary initial conditions for infinite Brownian particle systems within a large class and provides a construction of the particle systems themselves started from such initial conditions. Examples of particle systems falling into our framework include Brownian versions of TASEP-like processes such as the diffusive scaling limit of the $q$-TASEP process. In this context the spacings between consecutive particles form infinite-dimensional versions of the softly reflected Brownian motions recently introduced in the finite-dimensional setting by O'Connell and Ortmann and are of independent interest. The proof of the main result is based on intertwining relations satisfied by the particle systems involved which can be regarded as infinite-dimensional analogues of the suitably generalized Burke's Theorem.  
\end{abstract}

\maketitle


\section{Introduction}

Over the past decade Brownian particle systems have played a crucial role in a variety of problems in probability theory. Examples include (among others) the Dyson Brownian motion which is of great importance in the study of universal phenomena in random matrix theory (see the survey \cite{EY} and the references therein), as well as in obtaining the infinite-dimensional dynamics of eigenvalues arising in the bulk and at the edge of the spectrum of large random matrices (see \cite{Ts} for the former and \cite{Os1}, \cite{Os2}, \cite{Os3} for the latter); the stochastic version of the dyadic model which helps to understand anomalous dissipation of energy in the context of the Euler equation of fluid dynamics (see \cite{BFM1}, \cite{BFM2}, \cite{BFM3}); the systems of sticky Brownian particles giving rise to novel stochastic flows of kernels (see \cite{HW}); and rank-based and volatility-stabilized models of stochastic portfolio theory offering insights into the degree of stability of the capital distribution curve in financial markets (\cite{BFK}, \cite{PP}, \cite{PS}, \cite{IPS}). Despite these developments a unified approach to existence, uniqueness and (quasi-)stationarity in the case of infinitely many particles is missing. 

\medskip

The aim of this paper is to explicitly identify families of quasi-stationary measures for a large class of infinite Brownian particle systems, that is, initial conditions for which the joint distribution of spacings between consecutive particles does not change over time, and to give a construction of the particle systems started from such initial conditions. In addition, a martingale problem characterization is provided, yielding a powerful tool for establishing convergence to the infinite Brownian particle systems. The spacings between consecutive particles in such systems can be viewed as infinite-dimensional versions of the softly reflected Brownian motions recently introduced in the finite-dimensional setting by O'Connell and Ortmann in \cite{OO} and are of independent interest.

\medskip

Specifically, we consider systems of infinitely many particles on the real line whose positions $X_k$, $k\in\nn$ evolve according to the system of stochastic differential equations (SDEs)
\begin{equation}\label{main_sde}
\mathrm{d}X_k(t)=\mathrm{d}B_k(t)+\mu_k\,\mathrm{d}t
+\sum_{l=k}^{k+d-1} U'\big(X_l(t)-X_{l+1}(t)\big)\,r_{lk}\,\mathrm{d}t,\quad k\in\nn.
\end{equation}
Here $\{B_k:\,k\in\nn\}$ is an infinite system of one-dimensional Brownian motions with quadratic covariations $\langle B_k,B_l\rangle(t)=a_{kl}\,t$, $t\ge0$ where $a_{kl}$, $(k,l)\in\nn^2$ are real constants; $\mu_k$, $k\in\nn$ is a sequence of real numbers with $\mu_{k_0}=\mu_{k_0+1}=\cdots$ for some $k_0\in\nn$; $U:\,[0,\infty)\rightarrow\rr$ is a differentiable function with a locally bounded derivative; and $r_{lk}$, $(l,k)\in\nn^2$, $l\ge k$ are real constants. Note that the form of \eqref{main_sde} assumes that the interaction between particles is \textit{hierarchical}, that is, only particles with a higher index enter into the dynamics of any given particle, and of \textit{finite range}, that is, only finitely many particles appear in the dynamics of any given particle. Since a shift of a solution to \eqref{main_sde} gives another solution of \eqref{main_sde}, it may and will be assumed without loss of generality that $X_1(0)=0$. 

\medskip

Our main assumption is on the quantities $a_{kl}$, $(k,l)\in\nn^2$ and $r_{lk}$, $(l,k)\in\nn^2$, $l\ge k$ and should be viewed as the infinite-dimensional version of the \textit{skew-symmetry} condition in \cite{OO}. For notational convenience we use the convention $r_{lk}=0$ for $l\ge k+d$ and $l<k$. 

\begin{asmp}\label{OO_asmp}
The covariance matrix $A=(a_{kl})_{k,l=1}^\infty$ is such that all submatrices $(a_{kl})_{k,l=1}^K$, $K\in\nn$ are non-degenerate and 
\begin{eqnarray}
&& a_{kl}+a_{(k+1)(l+1)}-a_{(k+1)l}-a_{k(l+1)}=\frac{r_{lk}-r_{l(k+1)}}{2},\quad (k,l)\in\nn^2,\;l>k, \\
&& a_{1k}-a_{1(k+1)}=\frac{r_{k1}}{2},\quad k\in\nn.
\end{eqnarray}
We also assume the normalizations $a_{kk}+a_{(k+1)(k+1)}-a_{(k+1)k}-a_{k(k+1)}=1$, $k\in\nn$ and $r_{kk}=1$, $k\in\nn$.    
\end{asmp}

\begin{rmk}\label{OO_rmk}
Assumption \ref{OO_asmp} is satisfied if, for instance, $2A$ is the $\nn\times\nn$ identity matrix and $d=1$ (so that $r_{lk}=0$ for all $(l,k)\in\nn^2$, $l>k$). This simple setting includes Brownian versions of TASEP-like particle systems, and we refer to Section \ref{sec_examples} below for a detailed discussion of some examples of this type. 
\end{rmk}

\begin{rmk}
With $\tilde{A}:=\big(a_{kl}+a_{(k+1)(l+1)}-a_{(k+1)l}-a_{k(l+1)}\big)_{k,l=1}^\infty$ and $\tilde{R}:=\big(r_{lk}-r_{l(k+1)}\big)_{k,l=1}^\infty$ it is easy to see from Assumption \ref{OO_asmp} that the $\nn\times\nn$ matrix $2\tilde{A}-\tilde{R}$ is an upper triangular band matrix with band width $d$. In particular, there is a unique upper triangular matrix $\big(2\tilde{A}-\tilde{R}\big)^{-1}$ whose formal (left and right) multiplication with $2\tilde{A}-\tilde{R}$ gives the $\nn\times\nn$ identity matrix. Indeed, viewing the matrix $\big(2\tilde{A}-\tilde{R}\big)^{-1}$ as a left inverse first, one can fill in its entries from left to right and from top to bottom sequentially by solving the associated linear equations. Then, one checks that the resulting matrix is also a right inverse of $2\tilde{A}-\tilde{R}$. For future reference we set
\begin{equation}\label{nu_const}
(\nu_1,\nu_2,\ldots)^T=\big(2\tilde{A}-\tilde{R}\big)^{-1}\,(\mu_1-\mu_2,\mu_2-\mu_3,\ldots)^T
\end{equation}
where the superscript ``T'' denotes transpose.
\end{rmk}

In addition, we make mild assumptions on the interaction potential $U$.

\begin{asmp}\label{MP_asmp}
It holds ${\mathcal Z}_k:=\int_\rr \exp\big(2\,U(z)-2\,\nu_k\,z\big)\,\mathrm{d}z<\infty$, $k\in\nn$ and the probability measures $\frac{1}{{\mathcal Z}_k}\,\exp\big(2\,U(z)-2\,\nu_k\,z\big)\,\mathrm{d}z$, $k\in\nn$ have finite second moment and Fisher information:
\begin{equation}\label{Fisher}
\int_\rr \frac{1}{{\mathcal Z}_k}\,\exp\big(2\,U(z)-2\,\nu_k\,z\big)\,(2\,U'(z)-2\,\nu_k)^2\,\mathrm{d}z<\infty,
\quad k\in\nn\,.
\end{equation}
\end{asmp}

\begin{rmk}
As will become apparent below, the condition ${\mathcal Z}_k<\infty$, $k\in\nn$ ensures the existence of finite quasi-stationary measures for the process in \eqref{main_sde}. Moreover, it is easy to see that condition \eqref{Fisher} is equivalent to
\begin{equation}\label{Fisher'}
\int_\rr \exp\big(2\,U(z)-2\,\nu_k\,z\big)\,U'(z)^2\,\mathrm{d}z<\infty,\quad k\in\nn.
\end{equation}
\end{rmk}

\smallskip

The first main result of the paper reads as follows.

\begin{theorem}\label{main_thm1}
Suppose that Assumptions \ref{OO_asmp} and \ref{MP_asmp} are satisfied. Then there exists a unique Markovian weak solution of the infinite system of SDEs \eqref{main_sde} whose one-dimensional distributions are given by
\begin{equation}\label{1D_inf}
\begin{split}
&\prod_{k=1}^\infty \frac{1}{{\mathcal Z}_k}\,\exp\big(2\,U(x_k-x_{k+1})-2\,\nu_k\,(x_k-x_{k+1})\big)
\,\mathrm{d}(x_k-x_{k+1})\\
&\quad\quad\;\;\cdot\frac{1}{(2\pi a_{11} t)^{1/2}}\,\exp\bigg(-\frac{\big(x_1-\mu_1\,t-\sum_{l=1}^d r_{l1}\,\nu_l\,t\big)^2}{2 a_{11} t}\bigg)\,\mathrm{d}x_1,\quad t\ge0.
\end{split}
\end{equation} 
\end{theorem}

\smallskip

Picking a constant $\mu\in\rr$, subtracting $\mu d$ from all $\mu_k$, $k\in\nn$, and replacing $U$ in \eqref{main_sde} by $U_\mu(z):=U(z)+\mu z$ (clearly, this has no effect on \eqref{main_sde}) one obtains the following corollary.    
\begin{corollary}
Suppose that Assumption \ref{OO_asmp} holds. Then, for all $\mu\in\rr$ such that $U_\mu$ satisfies Assumption \ref{MP_asmp}, the measure 
\begin{equation}
\delta_0(\mathrm{d}x_1)\,\prod_{k=1}^\infty \frac{1}{{\mathcal Z}_{k,\mu}}
\,\exp\big(2\,U_\mu(x_k-x_{k+1})-2\,\nu_k\,(x_k-x_{k+1})\big)\,\mathrm{d}(x_k-x_{k+1})
\end{equation}
is quasi-stationary for the dynamics of \eqref{main_sde}. 
\end{corollary}

\smallskip

In the course of the proof of Theorem \ref{main_thm1} we establish the following result of independent interest. 

\begin{theorem}\label{main_thm2}
Suppose that Assumptions \ref{OO_asmp} and \ref{MP_asmp} are satisfied and consider a weak solution of \eqref{main_sde} as in Theorem \ref{main_thm1}. Then, for any fixed $K\in\nn$, the process $(X_k:\;1\le k\le K)$ constitutes the unique weak solution of the system of SDEs
\begin{equation}\label{main_sde_fin}
\mathrm{d}X_k(t)=\mathrm{d}W_k(t)+\mu_k\,\mathrm{d}t
+\sum_{l=k}^{K-1} U'\big(X_l(t)-X_{l+1}(t)\big)\,r_{lk}\,\mathrm{d}t
+\sum_{l=K}^{k+d-1} \nu_l\,r_{lk}\,\mathrm{d}t,\;\;\;1\le k\le K
\end{equation}
in its own filtration. Here $(W_k:\;1\le k\le K)$ is a Brownian motion with covariance matrix $(a_{kl})_{k,l=1}^K$. 
\end{theorem}

\smallskip

The statement of Theorem \ref{main_thm2} can be interpreted as an \textit{averaging principle}. Indeed, in the natural filtration of the Brownian particle system $(X_k:\;k\in\nn)$ the process $(X_k:\;1\le k\le K)$ simply satisfies the first $K$ SDEs in \eqref{main_sde}, whereas in the smaller filtration generated by $(X_k:\;1\le k\le K)$ (that is, after averaging out the interactions with the processes $X_l$, $l>K$) it solves \eqref{main_sde_fin}. That is, as a result of the averaging the interaction terms with the processes $X_l$, $l>K$ get replaced by suitable constants in the dynamics of the processes $X_k$, $1\le k\le K$.   

\medskip

Theorem \ref{main_thm2} provides a powerful tool for proving weak convergence to the solution of \eqref{main_sde} of Theorem \ref{main_thm1}. Indeed, in order to prove weak convergence of a sequence of processes $(\Xi^M_k:\;k\in\nn)$, $M\in\nn$ to the latter, it suffices to show weak convergence of each of the sequences $(\Xi^M_k:\;1\le k\le K)$, $M\in\nn$ to the solution of \eqref{main_sde_fin} with the same value of $K$. This, in turn, can be accomplished using standard martingale problem techniques (see e.g. \cite[Section 4.8]{EK}).   

\medskip

The proofs of Theorems \ref{main_thm1} and \ref{main_thm2} rely heavily on intertwinings of diffusion processes as discussed in great generality in \cite{PS2}. In fact, as the proofs show, for each $K\in\nn$, the diffusion of \eqref{main_sde_fin} is intertwined with the weak solution of \eqref{main_sde} described in Theorem \ref{main_thm1}. Moreover, any two diffusions of \eqref{main_sde_fin} with different values of $K$ are intertwined with each other as well. 

\medskip

The remainder of the paper is structured as follows. Section \ref{sec_ex} is devoted to the proof of the existence part in Theorem \ref{main_thm1}. The proof is based on the consistency of solutions to \eqref{main_sde_fin} with different values of $K$ which is the subject of Proposition \ref{prop_cons} below. Section \ref{sec_uniq} deals with the uniqueness part in Theorem \ref{main_thm1}. The proof of the latter is based on Theorem \ref{main_thm2} the proof of which, in turn, makes use of the entropy approach to time reversal developed by F\"ollmer and Wakolbinger in the series of papers \cite{Fo1}, \cite{Fo2}, \cite{FW} and is given in Section \ref{sec_uniq} as well. Finally, Section \ref{sec_examples} provides some examples from the area of integrable probability falling into the framework of \eqref{main_sde}.   

\section{Existence}\label{sec_ex}

A key ingredient in the proof of Theorem \ref{main_thm1} is the following proposition.

\begin{proposition}\label{prop_cons}
Suppose that Assumptions \ref{OO_asmp} and \ref{MP_asmp} are satisfied. Then:
\begin{enumerate}[(a)]
\item For each $K\in\nn$, there exists a unique weak solution of \eqref{main_sde_fin} started according to 
\begin{equation}\label{main_sde_fn_ic}
\delta_0(\mathrm{d}x_1)\,\prod_{k=1}^{K-1} \frac{1}{{\mathcal Z}_k}
\,\exp\big(2\,U(x_k-x_{k+1})-2\,\nu_k\,(x_k-x_{k+1})\big)\,\mathrm{d}(x_k-x_{k+1})
\end{equation} 
and, for $t>0$, its one-dimensional distributions are given by
\begin{equation}\label{1D_dist}
\frac{1}{(2\pi a_{11} t)^{1/2}}
\,\exp\bigg(-\frac{\big(x_1-\mu_1\,t-\sum_{l=1}^d r_{l1}\,\nu_l\,t\big)^2}{2 a_{11} t}\bigg)\,\mathrm{d}x_1
\,\prod_{k=1}^{K-1} \frac{1}{{\mathcal Z}_k}\,\exp\big(2\,U(y_k)-2\,\nu_k\,y_k\big)\,\mathrm{d}(x_k-x_{k+1}).
\end{equation}
\item The solution in part (a) has the property that, for any $1\le J<K$, the process $(X^{(K)}_k:\;1\le k\le J)$ is a weak solution of \eqref{main_sde_fin} with $K=J$. 
\end{enumerate}    
\end{proposition}

Proposition \ref{prop_cons} yields the existence part of Theorem \ref{main_thm1}. Indeed, it shows that the sequence of solutions $X^{(K)}$, $K\in\nn$ of \eqref{main_sde_fin} with initial conditions of the form \eqref{main_sde_fn_ic} is consistent. Hence, one can define a process $X=(X_k,\;k\in\nn)$ as the limit of $X^{(K)}$, $K\in\nn$ in the sense of the Kolmogorov Extension Theorem (see e.g. \cite[Theorem 5.16]{Ka}). Due to the fact that the $\sigma$-algebras in the natural filtration of $X$ are generated by cylinder events, it follows that the processes 
\begin{equation*}
X_k(t)-X_k(0)-\mu_k\,t-\int_0^t \sum_{l=k}^{k+d-1} U'\big(X_l(s)-X_{l+1}(s)\big)\,r_{lk}\,\mathrm{d}s,\quad k\in\nn
\end{equation*}
are continuous martingales in that filtration. Moreover, the quadratic covariations of the latter processes are given by $a_{kl}\,t$, $(k,l)\in\nn^2$ and, thus, $X$ is a weak solution of \eqref{main_sde_fin} by L\'evy's characterization of Brownian motion (see e.g. \cite[Theorem 3.3.16]{KS}). Lastly, $X$ has the Markov property, since, for any $t>0$, the $\sigma$-algebras generated by $X$ on $[0,t)$, at $t$, and on $(t,\infty)$ can be generated by cylinder events, and each $X^{(K)}$ has the Markov property (as a usual consequence of weak uniqueness for any initial condition given by \eqref{1D_dist}, see e.g. \cite[Theorem 5.4.20]{KS}).    

\medskip

Next, we give the proof of Proposition \ref{prop_cons}. 

\medskip

\noindent\textbf{Proof of Proposition \ref{prop_cons}. Step 1.} We fix a $K\in\nn$ throughout and start with the proof of part (a) in the case that $U'$ is bounded. The existence and uniqueness of a weak solution $X^{(K)}$ to \eqref{main_sde_fin} in this case is a well-known consequence of Girsanov's Theorem (see e.g. \cite[Section 5.3 B]{KS}). To complete the proof of part (a) in this case it therefore suffices to show that, for any $t>0$, the distribution of the random vector $\big(X^{(K)}_1(t),\,X^{(K)}_1(t)-X^{(K)}_2(t),\,\ldots,\,X^{(K)}_{K-1}(t)-X^{(K)}_K(t)\big)$ admits the density 
\begin{equation*}
\psi(t,x_1,y_1,\ldots,y_{K-1}):=\frac{1}{(2\pi a_{11} t)^{1/2}}
\,\exp\bigg(-\frac{\big(x_1-\mu_1\,t-\sum_{l=1}^d r_{l1}\,\nu_l\,t\big)^2}{2 a_{11} t}\bigg)
\,\prod_{k=1}^{K-1} \frac{1}{{\mathcal Z}_k}\,\exp\big(2\,U(y_k)-2\,\nu_k\,y_k\big).
\end{equation*}
To this end, we set 
\begin{eqnarray*}
&& \begin{split}
{\mathcal A}_{K-1}=
\frac{1}{2}\,\sum_{k,l=1}^{K-1} \big(a_{kl}+a_{(k+1)(l+1)}-a_{(k+1)l}-a_{k(l+1)}\big)
\,\frac{\partial^2}{\partial y_k\,\partial y_l}
+\sum_{k=1}^{K-1} \big(\mu_k-\mu_{k+1}\big)\,\frac{\partial}{\partial y_k} \\
+\sum_{k=1}^{K-1} \sum_{l=k}^{K-1} U'(y_l)\,\big(r_{lk}-r_{l(k+1)}\big)\,\frac{\partial}{\partial y_k} 
+\sum_{k=1}^{K-1} \Big(\sum_{l=K}^{k+d-1} \nu_l\,r_{lk}-\sum_{l=K}^{k+d} \nu_l\,r_{l(k+1)}\Big)
\,\frac{\partial}{\partial y_k},
\end{split} \\
&& \begin{split}
{\mathcal L}_K={\mathcal A}^{(K-1)}+\frac{1}{2}\,a_{11}\,\frac{\partial^2}{\partial x_1^2}
+\mu_1\,\frac{\partial}{\partial x_1}
+\sum_{l=1}^{K-1} U'(y_l)\,r_{l1}\,\frac{\partial}{\partial x_1}
+\sum_{l=K}^d \nu_l\,r_{lk}\,\frac{\partial}{\partial x_1} \\
+\sum_{k=1}^{K-1} (a_{1k}-a_{1(k+1)})\,\frac{\partial^2}{\partial x_1\partial y_k}
\end{split}
\end{eqnarray*}
and aim to verify the Kolmogorov forward equation $\frac{\partial\psi}{\partial t}=({\mathcal L}_K)^*\psi$. Here $({\mathcal L}_K)^*$ is the formal adjoint of ${\mathcal L}_K$, viewed as an operator acting on twice continuously differentiable functions.

\medskip

A direct computation gives 
\begin{equation*}
\frac{\partial\psi}{\partial t}=\bigg(-\frac{1}{2t}+\frac{x_1^2}{2a_{11}t^2}
-\frac{\big(\mu_1+\sum_{l=1}^d r_{l1}\,\nu_l\big)^2}{2a_{11}}\bigg)\,\psi.
\end{equation*}
Next, we note that 
\begin{equation*}
\begin{split}
({\mathcal L}_K)^*\psi=({\mathcal A}^{(K-1)})^*\psi
+\frac{1}{2}\,a_{11}\,\frac{\partial^2\psi}{\partial x_1^2}
-\mu_1\,\frac{\partial\psi}{\partial x_1}
-\sum_{l=1}^{K-1} U'(y_l)\,r_{l1}\,\frac{\partial\psi}{\partial x_1}
-\sum_{l=K}^d \nu_l\,r_{l1}\,\frac{\partial\psi}{\partial x_1} \\
+\sum_{k=1}^{K-1} (a_{1k}-a_{1(k+1)})\,\frac{\partial^2\psi}{\partial x_1\partial y_k}.
\end{split}
\end{equation*}
At this point, \cite[Theorem 3.2]{OO} in conjunction with Assumption \ref{OO_asmp} yield 
\begin{equation*}
({\mathcal A}^{(K-1)})^*\prod_{k=1}^{K-1} \frac{1}{{\mathcal Z}_k}\,\exp\big(2\,U(y_k)-2\,\nu_k\,y_k\big)=0,
\end{equation*}
so that $({\mathcal A}^{(K-1)})^*\psi=0$. Calculating the remaining terms directly we get
\begin{equation*}
\begin{split}
({\mathcal L}_K)^*\psi=&\bigg(-\frac{1}{2t}+\frac{\big(x_1-\mu_1\,t-\sum_{l=1}^d r_{l1}\,\nu_l\,t\big)^2}{2a_{11}t^2}\bigg)\,\psi
+\mu_1\,\frac{x_1-\mu_1\,t-\sum_{k=1}^d r_{k1}\,\nu_k\,t}{a_{11}t}\,\psi \\
&+\sum_{l=1}^{K-1} U'(y_l)\,r_{l1}\,\frac{x_1-\mu_1\,t-\sum_{k=1}^d r_{k1}\,\nu_k\,t}{a_{11}t}\,\psi 
+\sum_{l=K}^d \nu_l\,r_{l1}\,\frac{x_1-\mu_1\,t-\sum_{k=1}^d r_{k1}\,\nu_k\,t}{a_{11}t}\,\psi \\
&-\sum_{k=1}^{K-1} (a_{1k}-a_{1(k+1)})\,\big(2\,U'(y_k)-2\,\nu_k\big)
\,\frac{x_1-\mu_1\,t-\sum_{l=1}^d r_{l1}\,\nu_l\,t}{a_{11}t}\,\psi.
\end{split}
\end{equation*}
Using Assumption \ref{OO_asmp} this can be simplified further to
\begin{equation*}
\begin{split}
& \bigg(-\frac{1}{2t}+\frac{\big(x_1-\mu_1\,t-\sum_{l=1}^d r_{l1}\,\nu_l\,t\big)^2}{2a_{11}t^2}
+\Big(\mu_1+\sum_{k=1}^d r_{k1}\,\nu_k\Big)\,\frac{x_1-\mu_1\,t-\sum_{l=1}^d r_{l1}\,\nu_l\,t}{a_{11}t}\bigg)\,\psi \\
& = \bigg(-\frac{1}{2t}+\frac{x_1^2}{2a_{11}t^2}
-\frac{\big(\mu_1+\sum_{l=1}^d r_{l1}\,\nu_l\big)^2}{2a_{11}}\bigg)\,\psi.
\end{split}
\end{equation*}
Hence, $\frac{\partial\psi}{\partial t}=({\mathcal L}_K)^*\psi$ as desired. In view of the established weak uniqueness and \cite[Corollary 1.3, implication a) $\Rightarrow$ c)]{Ku} it follows that, for any $t>0$, the distribution of the process $\big(X^{(K)}_1,\,X^{(K)}_1-X^{(K)}_2,\ldots,X^{(K)}_{K-1}-X^{(K)}_K\big)$ at time $t$ has density $\psi(t,\cdot)$.  

\medskip

\noindent\textbf{Step 2.} We now turn to the proof of part (a) in the case that $U'$ is unbounded. To establish weak existence we consider a probability space which supports a $K$-dimensional Brownian motion $\beta$ with zero drift vector, covariance matrix $A^{(K)}:=(a_{kl})_{k,l=1}^K$ and initial condition of \eqref{main_sde_fn_ic}. Furthermore, we define the stopping times 
\begin{equation}\label{what_is_tau_n}
\tau_n:=\inf\Big\{s\ge0:\;\big|\beta_l(s)-\beta_{l+1}(s)\big|\ge n\;\;\text{for  some }\;1\le l\le K-1\Big\},\quad n\in\nn,
\end{equation}
fix a $t\in(0,\infty)$, and introduce the measures $\pp_n$, $n\in\nn$ given by their densities with respect to the underlying probability measure:
\begin{equation*}
\begin{split}
\frac{\mathrm{d}\pp_n}{\mathrm{d}\qq}=
\exp\bigg(\int_0^{\min(t,\tau_n)} \sum_{k,m=1}^K \Big(\mu_k+\sum_{l=k}^{K-1} U'\big(\beta_l(s)-\beta_{l+1}(s)\big)\,r_{lk}
+\sum_{l=K}^{k+d-1} \nu_l\,r_{lk}\Big)\,\big(A^{(K)}\big)^{-1}_{k,m}\,\mathrm{d}\beta_m(s) \qquad \\
-\frac{1}{2}\,\int_0^{\min(t,\tau_n)} \Big\|\Big(\mu_k+\sum_{l=k}^{K-1} U'\big(\beta_l(s)-\beta_{l+1}(s)\big)\,r_{lk}
+\sum_{l=K}^{k+d-1} \nu_l\,r_{lk}\Big)_{1\le k\le K}\Big\|^2_{(A^{(K)})^{-1}}\,\mathrm{d}s\bigg),\;\;\;n\in\nn.
\end{split}
\end{equation*}
Here $(A^{(K)})^{-1}$ is the inverse of $A^{(K)}$, and $\|\cdot\|_{(A^{(K)})^{-1}}$ stands for the Euclidean norm on $\rr^K$ associated with $(A^{(K)})^{-1}$. Since $U'$ is locally bounded, $\pp_n$, $n\in\nn$ are well-defined probability measures by the Novikov criterion (see e.g. \cite[Proposition 3.5.12]{KS}). 

\medskip

In addition, one computes 
\begin{equation*}
\begin{split}
& \ev^\qq\Big[\frac{\mathrm{d}\pp_n}{\mathrm{d}\qq}\,\log\frac{\mathrm{d}\pp_n}{\mathrm{d}\qq}\Big]
=\ev^{\pp_n}\Big[\log\frac{\mathrm{d}\pp_n}{\mathrm{d}\qq}\Big] \\
& =\ev^{\pp_n}\Big[\frac{1}{2}\,\int_0^{\min(t,\tau_n)} \Big\|\Big(\mu_k+\sum_{l=k}^{K-1} U'\big(\beta_l(s)-\beta_{l+1}(s)\big)\,r_{lk}
+\sum_{l=K}^{k+d-1} \nu_l\,r_{lk}\Big)_{1\le k\le K}\Big\|^2_{(A^{(K)})^{-1}}\,\mathrm{d}s\Big],\;\;n\in\nn
\end{split}
\end{equation*}
where the second equality is a consequence of Girsanov's Theorem. At this point, for every fixed $n\in\nn$, one can find a differentiable function $U_n$ which coincides with $U$ on $[-n,n]$ such that $U_n'$ is bounded, $\mathcal{Z}_{n,k}:=\int_\rr \exp\big(2\,U_n(z)-2\,\nu_k\,z\big)\,\mathrm{d}z<\infty$, $1\le k\le K-1$, and 
\begin{equation}\label{Fisher_control}
\int_\rr \frac{1}{\mathcal{Z}_{n,k}}\exp\big(2\,U_n(z)-2\,\nu_k\,z\big)\,U_n'(z)^2\,\mathrm{d}z
\le 1+\int_\rr \frac{1}{\mathcal{Z}_k}\exp\big(2\,U(z)-2\,\nu_k\,z\big)\,U'(z)^2\,\mathrm{d}z,\quad 1\le k\le K-1.
\end{equation}
Due to the weak uniqueness assertion established in Step 1, under $\pp_n$, the functional 
\begin{equation*}
\frac{1}{2}\,\int_0^{\min(t,\tau_n)} \Big\|\Big(\mu_k+\sum_{l=k}^{K-1} U'\big(\beta_l(s)-\beta_{l+1}(s)\big)\,r_{lk}
+\sum_{l=K}^{k+d-1} \nu_l\,r_{lk}\Big)_{1\le k\le K}\Big\|^2_{(A^{(K)})^{-1}}\,\mathrm{d}s
\end{equation*}
has the same law as the functional 
\begin{equation*}
\frac{1}{2}\,\int_0^{\min(t,\tau_n)} \Big\|\Big(\mu_k+\sum_{l=k}^{K-1} U_n'\big(X^{(K),n}_l(s)-X^{(K),n}_{l+1}(s)\big)\,r_{lk}+\sum_{l=K}^{k+d-1} \nu_l\,r_{lk}\Big)_{1\le k\le K}\Big\|^2_{(A^{(K)})^{-1}}\,\mathrm{d}s
\end{equation*}
where $X^{(K),n}$ is the weak solution of \eqref{main_sde_fin} with $U'$ replaced by $U_n'$ and initial condition
\begin{equation*}
\delta_0(\mathrm{d}x_1)\,\prod_{k=1}^{K-1} \frac{1}{{\mathcal Z}_{n,k}}
\,\exp\big(2\,U_n(x_k-x_{k+1})-2\,\nu_k\,(x_k-x_{k+1})\big)\,\mathrm{d}(x_k-x_{k+1}).
\end{equation*}
Relying on the formula for the one-dimensional distributions of $X^{(K),n}$ obtained in Step 1 one can now bound the expectation of the latter functional from above by a constant independent of $n$ using  
\begin{equation*}
\ev\bigg[\int_0^t U_n'\big(X^{(K),n}_l(s)-X^{(K),n}_{l+1}(s)\big)^2\,\mathrm{d}s\bigg]
=t\,\int_\rr \frac{1}{\mathcal{Z}_{n,l}}\exp\big(2\,U_n(z)-2\,\nu_l\,z\big)\,U_n'(z)^2\,\mathrm{d}z,\quad 1\le l\le K-1
\end{equation*}
and \eqref{Fisher_control}. 

\medskip

All in all, it follows that the quantities $\ev^\qq\big[\frac{\mathrm{d}\pp_n}{\mathrm{d}\qq}\,\log\frac{\mathrm{d}\pp_n}{\mathrm{d}\qq}\big]$, $n\in\nn$ can be bounded above uniformly in $n$, so that $\frac{\mathrm{d}\pp_n}{\mathrm{d}\qq}$, $n\in\nn$ is a uniformly integrable martingale under $\qq$ (the martingale property can be deduced directly from It\^o's formula). Thus, Doob's Martingale Convergence Theorem yields the existence of the limit $\frac{\mathrm{d}\pp}{\mathrm{d}\qq}:=\lim_{n\to\infty} \frac{\mathrm{d}\pp_n}{\mathrm{d}\qq}$ in the $L^1$ sense, and we can let $\pp$ be the probability measure associated with the density $\frac{\mathrm{d}\pp}{\mathrm{d}\qq}$. Moreover, the $L^1$ convergence implies that, for every twice continuously differentiable function $f:\,\rr^K\to\rr$ which is bounded together with all its first and second order partial derivatives, $n\in\nn$, $s_1,s_2\in[0,t]$, and bounded random variable $\Lambda$ measurable with respect to the $\sigma$-algebra generated by $\beta(\min(s,\tau_n))$, $s\in[0,s_1]$, one has
\begin{equation*}
\ev^\qq\bigg[\frac{\mathrm{d}\pp}{\mathrm{d}\qq}\bigg(f\big(\beta\big(\min(s_2,\tau_n)\big)\big)
-f\big(\beta\big(\min(s_1,\tau_n)\big)\big)
-\int_{\min(s_1,\tau_n)}^{\min(s_2,\tau_n)} \big(\widetilde{\mathcal L}_K f\big)\big(\beta\big(\min(s,\tau_n\big)\big)\big)
\,\mathrm{d}s\bigg)\,\Lambda\bigg]=0
\end{equation*} 
where $\widetilde{\mathcal L}_K$ is the generator corresponding to the SDE \eqref{main_sde_fin}. In other words, under $\pp$, the process $\beta$ solves the local martingale problem associated with \eqref{main_sde_fin} on $[0,t]$ and is, hence, a weak solution of \eqref{main_sde_fin} on $[0,t]$ (see e.g. \cite[Corollary 5.4.8]{KS}). The weak existence on the whole time interval $[0,\infty)$ follows from the Kolmogorov Extension Theorem and the weak uniqueness on every time interval $[0,t]$ established next. 

\medskip

In order to show weak uniqueness on each time interval $[0,t]$, we pick a weak solution $X^{(K)}$ of \eqref{main_sde_fin} on $[0,t]$, redefine the stopping times $\tau_n$, $n\in\nn$ by replacing $\beta_l(s)-\beta_{l+1}(s)$ with $X^{(K)}_l(s)-X^{(K)}_{l+1}(s)$ in \eqref{what_is_tau_n} and only considering $s\in[0,t]$ there, and note that 
\begin{equation}\label{L2controlwu}
\begin{split}
& \ev^\pp\bigg[\int_0^t \Big\|\Big(\mu_k+\sum_{l=k}^{K-1} U'\big(X^{(K)}_l(s)-X^{(K)}_{l+1}(s)\big)\,r_{lk}
+\sum_{l=K}^{k+d-1} \nu_l\,r_{lk}\Big)_{1\le k\le K}\Big\|^2_{(A^{(K)})^{-1}}\,\mathrm{d}s\bigg] \\
& =\lim_{n\to\infty} \,\ev^\pp\bigg[\int_0^{\tau_n} \Big\|\Big(\mu_k+\sum_{l=k}^{K-1} U'\big(X^{(K)}_l(s)-X^{(K)}_{l+1}(s)\big)\,r_{lk} +\sum_{l=K}^{k+d-1} \nu_l\,r_{lk}\Big)_{1\le k\le K}\Big\|^2_{(A^{(K)})^{-1}}\,\mathrm{d}s\bigg]
\end{split}
\end{equation}
by the Monotone Convergence Theorem. Moreover, the latter expectations can be bounded above uniformly in $n$ as before, by replacing $U$ with $U_n$ and $X^{(K)}$ with $X^{(K),n}$, changing the interval of integration from $[0,\tau_n]$ to $[0,t]$, and using the formula for the one-dimensional distributions of $X^{(K),n}$ obtained in Step 1. It follows that the expectation on the left-hand side of \eqref{L2controlwu} is finite. 

\medskip

At this point, \cite[Theorem 7.5]{LS} shows that the law of $X^{(K)}$ on $[0,t]$ is absolutely continuous with respect to the law of the previously introduced Brownian motion $\beta$ on $[0,t]$. Finally, the formula for the associated change of measure in \cite[Theorem 7.6, equation (7.29)]{LS} demonstrates that the law of $X^{(K)}$ is uniquely determined.  
 
\medskip

\noindent\textbf{Step 3.} To prove part (b) one only needs to consider the case $J=K-1$, since for all other values of $1\le J<K$ one can apply the statement for $J=K-1$ repeatedly. Set
\begin{eqnarray*}
&&Y^{(K)}=\big(X^{(K)}_1,\,X^{(K)}_1-X^{(K)}_2,\,\ldots,\,X^{(K)}_{K-1}-X^{(K)}_K\big), \\
&&Y^{(K,K-1)}=\big(X^{(K)}_1,\,X^{(K)}_1-X^{(K)}_2,\,\ldots,\,X^{(K)}_{K-2}-X^{(K)}_{K-1}\big),
\end{eqnarray*}
and let ${\mathcal F}(t)$, $t\ge0$ and ${\mathcal G}(t)$, $t\ge0$ be the filtrations generated by $Y^{(K)}$ and $Y^{(K,K-1)}$, respectively. 

\medskip

It suffices to show that the process 
\begin{equation}\label{whatisMf}
M^f(t)=f\big(Y^{(K,K-1)}(t)\big)-f\big(Y^{(K,K-1)}(0)\big)
-\int_0^t ({\mathcal L}_{K-1}f)\big(Y^{(K,K-1)}(s)\big)\,\mathrm{d}s,\quad t\ge0
\end{equation}
is a martingale with respect to ${\mathcal G}(t)$, $t\ge0$ for any twice continuously differentiable function $f:\,\rr^{K-1}\to\rr$ which is bounded together with all its first and second order partial derivatives. To this end, we employ the martingale property of the process
\begin{equation*}
f\big(Y^{(K,K-1)}(t)\big)-f\big(Y^{(K,K-1)}(0)\big)
-\int_0^t ({\mathcal L}_Kf)\big(Y^{(K)}(s)\big)\,\mathrm{d}s,\quad t\ge0
\end{equation*}
with respect to ${\mathcal F}(t)$, $t\ge0$ (due to It\^o's formula) to compute
\begin{eqnarray*}
\ev\big[f\big(Y^{(K,K-1)}(t)\big)\big|{\mathcal G}(s)\big]
&=& \ev\big[\ev\big[f\big(Y^{(K,K-1)}(t)\big)\big|{\mathcal F}(s)\big]\big|{\mathcal G}(s)\big] \\
&=& f\big(Y^{(K,K-1)}(s)\big)
+\ev\bigg[\ev\bigg[\int_s^t ({\mathcal L}_Kf)\big(Y^{(K)}(u)\big)\,\mathrm{d}u\bigg|{\mathcal F}(s)\bigg]
\bigg|{\mathcal G}(s)\bigg] \\
&=& f\big(Y^{(K,K-1)}(s)\big)+\ev\bigg[\int_s^t ({\mathcal L}_{K-1}f)\big(Y^{(K,K-1)}(u)\big)\,\mathrm{d}u
\bigg|{\mathcal G}(s)\bigg] \\
&&\quad\quad\quad\quad\quad\quad\quad\;
+\ev\bigg[\int_s^t \big(({\mathcal L}_K-{\mathcal L}_{K-1})f\big)\big(Y^{(K)}(u)\big)\,\mathrm{d}u
\bigg|{\mathcal G}(s)\bigg]
\end{eqnarray*}
for any $t\ge s\ge0$. Therefore we need to show that the last summand in the latter expression vanishes. Moreover, due to Fubini's Theorem, it is sufficient to prove
\begin{equation}\label{prop_claim1}
\ev\Big[\big(({\mathcal L}_K-{\mathcal L}_{K-1})f\big)\big(Y^{(K)}(u)\big)\Big|{\mathcal G}(s)\Big]=0,\quad 
u\ge s\ge 0. 
\end{equation}

\smallskip

We explain first how \eqref{prop_claim1} can be obtained from the following claim the proof of which we defer to Step 4 below.  

\medskip

\noindent\textbf{Claim 1.} For any $u\ge 0$, the $\sigma$-algebra ${\mathcal G}(u)$ and the random variable $Y^{(K)}_K(u)$ are independent. 

\medskip

Given the claim we find
\begin{eqnarray*}
\ev\Big[\big(({\mathcal L}_K-{\mathcal L}_{K-1})f\big)\big(Y^{(K)}(u)\big)\Big|{\mathcal G}(s)\Big] 
&=&\ev\Big[\ev\Big[\big(({\mathcal L}_K-{\mathcal L}_{K-1})f\big)\big(Y^{(K)}(u)\big)
\Big|{\mathcal G}(u)\Big]\Big|{\mathcal G}(s)\Big] \\
&=&\ev\Big[\ev\Big[\big(({\mathcal L}_K-{\mathcal L}_{K-1})f\big)\big(Y^{(K)}(u)\big)
\Big|Y^{(K,K-1)}(u)\Big]\Big|{\mathcal G}(s)\Big]. 
\end{eqnarray*}
Hence, to obtain \eqref{prop_claim1} it suffices to prove
\begin{equation}\label{prop_claim2}
\ev\Big[\big(({\mathcal L}_K-{\mathcal L}_{K-1})f\big)\big(Y^{(K)}(u)\big)\Big|Y^{(K,K-1)}(u)\Big]=0,
\quad u\ge0.
\end{equation}
At this point, we use formula \eqref{1D_dist} to write
\begin{eqnarray*}
&&\ev\Big[\big(({\mathcal L}_K-{\mathcal L}_{K-1})f\big)\big(Y^{(K)}(u)\big)\Big|Y^{(K,K-1)}(u)\Big] \\
&&=\int_\rr \bigg(\frac{\partial f}{\partial x_1}\big(Y^{(K,K-1)}(u)\big)\,r_{(K-1)1}
+\sum_{l=1}^{K-2} \frac{\partial f}{\partial y_l}\big(Y^{(K,K-1)}(u)\big)
\,\big(r_{(K-1)l}-r_{(K-1)(l+1)}\big)\bigg) \\
&&\qquad\qquad\qquad\qquad\qquad\qquad\qquad\quad\quad\quad\;\;
\cdot\big(U'(z)-\nu_K\big)\,\frac{1}{{\mathcal Z}_K}\,\exp(2\,U(z)-2\,\nu_K\,z)\,\mathrm{d}z.
\end{eqnarray*} 
Since the latter integrand is a derivative in $z$ of a function vanishing at positive and negative infinity, its integral equals to zero. All in all, it follows that the process $M^f$ of \eqref{whatisMf} is a martingale with respect to ${\mathcal G}(t)$, $t\ge0$ for any function $f$ as described in the paragraph after \eqref{whatisMf}.  

\medskip

\noindent\textbf{Step 4.} To finish the proof it remains to establish Claim 1. To this end, it suffices to consider the time-reversed process $\widehat{Y}^{(K,K-1)}(s):=Y^{(K,K-1)}(u-s)$, $s\in[0,u]$ in the filtration $\widehat{\mathcal{F}}(s)$, $s\in[0,u]$ generated by $Y^{(K)}(u-s)$, $s\in[0,u]$ and to verify that its dynamics does not depend on $Y^{(K)}_K(u)$. 

\medskip

Recall from Steps 1 and 2 that there exists an equivalent probability measure $\qq$ on the underlying probability space under which $Y^{(K,K-1)}(s)$, $s\in[0,u]$ is a Brownian motion with zero drift vector, covariance matrix $\Theta^{(K-1)}=(\theta_{kl})_{k,l=1}^{K-1}$ with entries
\begin{equation*}
\theta_{kl}:=
\begin{cases}
a_{kl} &\text{if}\;\;\;k=l=1 \\
a_{k(l-1)}-a_{kl} &\text{if}\;\;\;k=1,\,l>1 \\
a_{(k-1)l}-a_{kl} &\text{if}\;\;\;k>1,\,l=1 \\
a_{(k-1)(l-1)}+a_{kl}-a_{k(l-1)}-a_{(k-1)l} &\text{if}\;\;\;k>1,\,l>1
\end{cases}
\end{equation*}
and initial condition
\begin{equation}\label{term_cond}
\delta_0(\mathrm{d}x_1)\,\prod_{k=1}^{K-2} \frac{1}{{\mathcal Z}_k}
\,\exp\big(2\,U(y_k)-2\,\nu_k\,y_k\big)\,\mathrm{d}y_k
\end{equation}
in the filtration ${\mathcal F}(s)$, $s\in[0,u]$. Consequently, under $\qq$ and in the filtration $\widehat{\mathcal{F}}(s)$, $s\in[0,u]$ the time-reversed process $\widehat{Y}^{(K,K-1)}(s)$, $s\in[0,u]$ is a Brownian motion with zero drift vector and covariance matrix $\Theta^{(K-1)}$ which is conditioned on having the distribution \eqref{term_cond} at time $u$. More specifically, $\widehat{Y}^{(K,K-1)}(s)$, $s\in[0,u]$ solves the SDE
\begin{equation}\label{gamma2_eq}
\begin{split}
\mathrm{d}\widehat{Y}^{(K,K-1)}(s)=\mathrm{d}\widehat{w}(s)+\Theta^{(K-1)}\,
\bigg(\nabla_y
\log\int_{\rr^{K-1}} \exp\bigg(-\frac{(y-z)^T\big(\Theta^{(K-1)}\big)^{-1}(y-z)}{2(u-s)}\bigg) \;\;\;\qquad\qquad\qquad\\
\,\delta_0(\mathrm{d}z_1)\,\prod_{l=2}^{K-1} \exp\big(2\,U(z_l)-2\,\nu_l\,z_l\big)\,\mathrm{d}z_l\bigg)
\bigg|_{y=\widehat{Y}^{(K,K-1)}(s)}\,\mathrm{d}s \\
=:\mathrm{d}\widehat{w}(s)+\gamma^{(2)}(s,\widehat{Y}^{(K,K-1)}(s))\,\mathrm{d}s
\end{split}
\end{equation}
where $\widehat{w}$ is a $(K-1)$-dimensional Brownian motion with zero drift vector and covariance matrix $\Theta^{(K-1)}$ under $\qq$ and in the filtration $\widehat{\mathcal{F}}(s)$, $s\in[0,u]$ (take e.g. both $X$ and $B$ to be $Y^{(K,K-1)}$ in \cite[Theorem 1]{Me}). 

\medskip

Note that the drift function $\gamma^{(2)}(s,y)$ can be rewritten as
\begin{equation*}
\begin{split}
-\frac{y}{u-s}+\Theta^{(K-1)}\,\nabla_y\log\int_{\rr^{K-1}} \exp\bigg(\frac{y^T\big(\Theta^{(K-1)}\big)^{-1}z}{u-s}\bigg)
\,\exp\bigg(-\frac{z^T\big(\Theta^{(K-1)}\big)^{-1}z}{2(u-s)}\bigg) \\
\,\delta_0(\mathrm{d}z_1)\,\prod_{l=2}^{K-1} \exp\big(2\,U(z_l)-2\,\nu_l\,z_l\big)\,\mathrm{d}z_l
\end{split}
\end{equation*}
and that the second summand in the latter expression is given by a linear transformation of the gradient of the logarithmic moment generating function of a finite measure with sub-Gaussian tails. Hence, $\gamma^{(2)}(s,y)$ is of uniform sub-linear growth in $y$ on compact subsets of $s\in[0,u)$. It now follows from \cite[Corollary 3.5.16]{KS} that the law of $\widehat{Y}^{(K,K-1)}(s)$, $s\in[0,u)$ under $\qq$ in the filtration $\widehat{\mathcal{F}}(s)$, $s\in[0,u)$ is locally absolutely continuous with respect to the law of $\widehat{w}$ under $\qq$ in the filtration $\widehat{\mathcal{F}}(s)$, $s\in[0,u)$ with the corresponding densities being given by
\begin{equation}\label{Gir1}
\begin{split}
\exp\bigg(\int_0^s \gamma^{(2)}(q,\widehat{Y}^{(K,K-1)}(q))^T\big(\Theta^{(K-1)}\big)^{-1}\mathrm{d}\hat{w}(q)
-\frac{1}{2}\int_0^s \big\|\gamma^{(2)}(q,\widehat{Y}^{(K,K-1)}(q))\big\|^2_{(\Theta^{(K-1)})^{-1}}\,\mathrm{d}q\bigg),\\
s\in[0,u).
\end{split}
\end{equation}

\smallskip

On the other hand, one can employ the key observation made in \cite[proof of Lemma 3.1]{Fo2}: since the law of the process $Y^{(K,K-1)}(s)$, $s\in[0,u]$ under the original probability measure $\pp$ has finite relative entropy with respect to the law of the Brownian motion with the same initial condition, zero drift vector and covariance matrix $\Theta^{(K-1)}$ (see the estimates on the quantity of \eqref{L2controlwu} above), the same is true for the laws of their time reversals. In other words, the law of the process $\widehat{Y}^{(K,K-1)}(s)$, $s\in[0,u]$ under $\pp$ must have finite relative entropy with respect to the law of the same process under $\qq$. In particular, we conclude that the law of the process $\widehat{Y}^{(K,K-1)}(s)$, $s\in[0,u)$ under $\pp$ in the filtration $\widehat{\mathcal{F}}(s)$, $s\in[0,u)$ is locally absolutely continuous with respect to the law of $\widehat{w}$ under $\qq$ in the filtration $\widehat{\mathcal{F}}(s)$, $s\in[0,u)$. By \cite[Theorem 7.11]{LS} it follows that the process $\widehat{Y}^{(K,K-1)}(s)$, $s\in[0,u)$ satisfies a SDE of the form
\begin{equation}\label{gamma_eq}
\mathrm{d}\widehat{Y}^{(K,K-1)}(s)=\gamma(s)\,\mathrm{d}s+\mathrm{d}\widehat{b}(s),\quad 0\le s<u
\end{equation}
under $\pp$ and in the filtration $\widehat{\mathcal{F}}(s)$, $s\in[0,u)$, and \cite[Theorems 7.5, 7.6]{LS} show that the densities associated with the latter absolute continuity relation are given by
\begin{equation}\label{Gir2}
\exp\bigg(\int_0^s \gamma(q)^T\big(\Theta^{(K-1)}\big)^{-1}\mathrm{d}\widehat{Y}^{(K,K-1)}(q)
-\frac{1}{2}\int_0^s \|\gamma(q)\|^2_{(\Theta^{(K-1)})^{-1}}\,\mathrm{d}q\bigg), \quad s\in[0,u).
\end{equation}
Here $\gamma$ and $\widehat{b}$ are a non-anticipative process and a Brownian motion with zero drift vector and covariance matrix $\Theta^{(K-1)}$, respectively, in the filtration $\widehat{\mathcal{F}}(s)$, $s\in[0,u)$. 

\medskip

To identify the process $\gamma$ (and, thus, to complete the proof of Claim 1) one can use the decomposition $\gamma(s)=:\gamma^{(1)}(s)+\gamma^{(2)}(s,\widehat{Y}^{(K,K-1)}(s))$, $s\in[0,u)$ and combine \eqref{Gir1}, \eqref{Gir2} and Girsanov's Theorem in the form of \cite[Theorem 3.5.1]{KS} to conclude that the densities associated with the change of measure from \eqref{gamma2_eq} to \eqref{gamma_eq} are given by
\begin{equation*}
\exp\bigg(\int_0^s \gamma^{(1)}(q)\big(\Theta^{(K-1)}\big)^{-1}\mathrm{d}\widehat{b}(q)
+\frac{1}{2}\int_0^s \big\|\gamma^{(1)}(q)\big\|^2_{(\Theta^{(K-1)})^{-1}}\,\mathrm{d}q\bigg),\quad s\in[0,u).
\end{equation*}
Hence, the finiteness of the corresponding relative entropy together with the simple localization argument in \cite[proof of Lemma 2.6]{Fo1} imply
\begin{equation*}
\ev^\pp\bigg[\int_0^s \big\|\gamma^{(1)}(q)\big\|^2_{(\Theta^{(K-1)})^{-1}}\,\mathrm{d}q\bigg]<\infty,\quad s\in[0,u)
\end{equation*}
and, thus,
\begin{equation*}
\ev^\pp\bigg[\int_0^s \big\|\gamma^{(1)}(q)\big\|^2\,\mathrm{d}q\bigg]<\infty,\quad s\in[0,u)
\end{equation*}
where $\|\cdot\|$ is the standard Euclidean norm on $\rr^{K-1}$. This estimate can be improved further to
\begin{equation}\label{a_priori}
\ev^\pp\bigg[\int_0^s \big\|\gamma(q)\big\|^2\,\mathrm{d}q\bigg]<\infty,\quad s\in[0,u)
\end{equation}
by using the decomposition $\gamma(q)=\gamma^{(1)}(q)+\gamma^{(2)}(q,\widehat{Y}^{(K,K-1)}(q))$, recalling the uniform sub-linear growth of $\gamma^{(2)}(q,\widehat{Y}^{(K,K-1)}(q))$ in $\widehat{Y}^{(K,K-1)}(q)$ on compact subsets of $q\in[0,u)$, and combining part (a) of the proposition with the finite second moment condition in Assumption \ref{MP_asmp}. 
 
\medskip

In view of the estimate \eqref{a_priori} one can now apply \cite[Proposition 2.5]{Fo2} (more precisely its straightforward extension to higher dimensions) to obtain the representation
\begin{equation*}
\gamma(s)=\lim_{\epsilon\downarrow0} \;
\epsilon^{-1}\;\ev^\pp\big[\widehat{Y}^{(K,K-1)}(s+\epsilon)-\widehat{Y}^{(K,K-1)}(s)\big|\widehat{\mathcal F}(s)\big]
=\lim_{\epsilon\downarrow0} \;\epsilon^{-1}\;\int_{s-\epsilon}^s \gamma(q)\,\mathrm{d}q
\end{equation*}
for Lebesgue almost every $s\in[0,u)$ where the limits should be understood in the $L^2$ sense. This representation can be simplified further to
\begin{equation}\label{gamma_char}
\gamma(s)=\lim_{\epsilon\downarrow0} \;
\epsilon^{-1}\;\ev^\pp\big[\widehat{Y}^{(K,K-1)}(s+\epsilon)-\widehat{Y}^{(K,K-1)}(s)\big|Y^{(K)}(u-s)\big]
=\lim_{\epsilon\downarrow0} \;\epsilon^{-1}\;\int_{s-\epsilon}^s \gamma(q)\,\mathrm{d}q
\end{equation}
by noting that $Y^{(K)}$ (and, hence, also $Y^{(K)}(u-s)$, $s\in[0,u)$) is a Markov process. Indeed, the proof of weak uniqueness in part (a) of the proposition can be repeated word for word for an initial condition as in \eqref{1D_dist}, and this yields the Markov property of $Y^{(K)}$ by a standard argument (see e.g. \cite[Theorem 5.4.20]{KS}).

\medskip

Now, for any twice continuously differentiable function $f:\,\rr^K\to\rr$ with compact support one has, on the one hand,
\begin{equation}\label{gamma1}
\lim_{\epsilon\downarrow0}\;\epsilon^{-1}\;\ev^\pp\Big[\big(\widehat{Y}^{(K,K-1)}(s+\epsilon)
-\widehat{Y}^{(K,K-1)}(s)\big) f\big(Y^{(K)}(u-s)\big)\Big]=\ev^\pp\big[\gamma(s)\,f\big(Y^{(K)}(u-s)\big)\big]
\end{equation}
for Lebesgue almost every $s\in[0,u)$ by virtue of \eqref{gamma_char}. On the other hand, writing $D^{(K,K-1)}$ for the drift coefficients of $Y^{(K,K-1)}$ under $\pp$ and in the filtration ${\mathcal F}(t)$, $t\ge0$, one computes
\begin{equation*}
\begin{split}
&\lim_{\epsilon\downarrow0}\;\epsilon^{-1}\;\ev^\pp\Big[\big(\widehat{Y}^{(K,K-1)}(s+\epsilon)
-\widehat{Y}^{(K,K-1)}(s)\big) f\big(Y^{(K)}(u-s)\big)\Big] \\
&=-\lim_{\epsilon\downarrow0}\;\epsilon^{-1}\;\ev^\pp\Big[Y^{(K,K-1)}(u-s)\,f\big(Y^{(K)}(u-s)\big)
-Y^{(K,K-1)}(u-s-\epsilon)\,f\big(Y^{(K)}(u-s-\epsilon)\big)\Big] \\
&\quad+\lim_{\epsilon\downarrow0}\;\epsilon^{-1}\;\ev^\pp\Big[Y^{(K,K-1)}(u-s-\epsilon)
\Big(f\big(Y^{(K)}(u-s)\big)-f\big(Y^{(K)}(u-s-\epsilon)\big)\Big)\Big] \\
&=-\ev^\pp\Big[Y^{(K,K-1)}(u-s)\,\big({\mathcal L}_K f\big)\big(Y^{(K)}(u-s)\big)\Big]
-\ev^\pp\Big[f\big(Y^{(K)}(u-s)\,D^{(K,K-1)}(u-s)\Big] \\
&\quad-\ev^\pp\bigg[\bigg(\frac{\partial f}{\partial x_1}\big(Y^{(K)}(u-s)\big)\,\Theta^{(K)}_{1k}
+\sum_{l=1}^{K-1} \frac{\partial f}{\partial y_l}\big(Y^{(K)}(u-s)\big)\,\Theta^{(K)}_{(l+1)k}\bigg)_{1\le k\le K-1}\bigg] \\ &\quad+\ev^\pp\Big[Y^{(K,K-1)}(u-s)\,\big({\mathcal L}_K f\big)\big(Y^{(K)}(u-s)\big)\Big] 
\end{split}
\end{equation*}
for Lebesgue almost every $s\in[0,u)$ using It\^o's formula and forward in time versions of \eqref{gamma_char}. By cancelling out the two identical expectations and integrating by parts one can simplify the latter expression to
\begin{equation}\label{gamma2}
\begin{split}
\ev^\pp\bigg[f\big(Y^{(K)}(u-s)\bigg(-D^{(K,K-1)}(u-s)
+\bigg(-\frac{Y^{(K)}_1(u-s)-\big(\mu_1+\sum_{l=1}^d r_{l1}\nu_l\big)(u-s)}{a_{11}(u-s)}\,\Theta^{(K)}_{1k} \quad\quad\quad\quad\;\;\;\;\,\\
+\sum_{l=2}^K \Big(2\,U'\big(Y^{(K)}_l(u-s)\big)-2\,\nu_l\Big)\,\Theta^{(K)}_{lk}\bigg)_{1\le k\le K-1}\bigg)\bigg].
\end{split}
\end{equation}
Since $\gamma(s)$ is measurable with respect to the $\sigma$-algebra generated by $Y^{(K)}(u-s)$, a comparison of \eqref{gamma1} and \eqref{gamma2} allows to identify $\gamma(s)$ for Lebesgue almost every $s\in[0,u)$. Moreover, Assumption \ref{OO_asmp} implies that, for such $s$, the drift coefficient $\gamma(s)$ is a deterministic function of $Y^{(K,K-1)}(u-s)=\widehat{Y}^{(K,K-1)}(s)$ only. Finally, using this and repeating the proof of weak uniqueness in Step 2 mutatis mutandis for the SDE resulting from \eqref{gamma_eq}, one concludes that the dynamics of the process $\widehat{Y}^{(K,K-1)}(s)$, $s\in[0,u]$ in the filtration $\widehat{\mathcal F}(s)$, $s\in[0,u]$ does not depend on $Y^{(K)}_K(u)$. Claim 1 readily follows. \ep 

\section{Uniqueness}\label{sec_uniq}

We first explain how Theorem \ref{main_thm2} gives the uniqueness part of Theorem \ref{main_thm1} and, thus, completes its proof. To this end, we consider two weak solutions $X$ and $\tilde{X}$ of \eqref{main_sde} as described in Theorem \ref{main_thm1}. By Theorem \ref{main_thm2}, for any $K\in\nn$, the processes $(X_k:\;1\le k\le K)$ and $(\tilde{X}_k:\;1\le k\le K)$ are both weak solutions to \eqref{main_sde_fin}. Therefore, by Proposition \ref{prop_cons} (a) their laws must be the same. It follows that $X$ and $\tilde{X}$ have the same law.

\medskip

Next, we prove Theorem \ref{main_thm2}. The proof is similar to the proof of Proposition \ref{prop_cons}, but requires certain modifications and is therefore presented in full detail for the convenience of the reader. 

\medskip

\noindent\textbf{Proof of Theorem \ref{main_thm2}. Step 1.} We start by fixing a $K\in\nn$, putting 
\begin{equation*}
Y^{(\infty)}:=\big(X_1,\,X_1-X_2,\,X_2-X_3,\,\ldots\big)\quad\text{and}\quad
Y^{(\infty,K)}:=\big(X_1,\,X_1-X_2,\,\ldots,\,X_{K-1}-X_K\big),
\end{equation*}
and letting ${\mathcal F}^{(\infty)}(t)$, $t\ge0$ and ${\mathcal F}^{(\infty,K)}(t)$, $t\ge0$ be the filtrations generated by $Y^{(\infty)}$ and $Y^{(\infty,K)}$, respectively. Moreover, we formally define the operator
\begin{eqnarray*}
{\mathcal L}_\infty&=&
\frac{1}{2}\,a_{11}\,\frac{\partial^2}{\partial x_1^2}
+\sum_{k=1}^\infty (a_{1k}-a_{1(k+1)})\,\frac{\partial^2}{\partial x_1\partial y_k}
+\mu_1\,\frac{\partial}{\partial x_1}
+\sum_{l=1}^\infty U'(y_l)\,r_{l1}\,\frac{\partial}{\partial x_1} \\
&&+\frac{1}{2}\,\sum_{k,l=1}^\infty \big(a_{kl}+a_{(k+1)(l+1)}-a_{(k+1)l}-a_{k(l+1)}\big)
\,\frac{\partial^2}{\partial y_k\,\partial y_l}
+\sum_{k=1}^\infty \big(\mu_k-\mu_{k+1}\big)\,\frac{\partial}{\partial y_k} \\
&&+\sum_{k=1}^\infty \sum_{l=k}^\infty U'(y_l)\,\big(r_{lk}-r_{l(k+1)}\big)\,\frac{\partial}{\partial y_k} 
\end{eqnarray*}
which should be understood as acting on functions $f:\,\rr^\nn\rightarrow\rr$ depending only on finitely many coordinates of $(x_1,\,y_1,\,y_2,\,\ldots)\in\rr^\nn$ and being twice continuously differentiable with respect to these coordinates. 

\medskip

Recalling the notation ${\mathcal L}_K$ introduced in Step 1 of the proof of Proposition \ref{prop_cons} we note that, in order to establish the theorem, it suffices to show that the process 
\begin{equation}\label{whatisMf'}
f\big(Y^{(\infty,K)}(t)\big)-f\big(Y^{(\infty,K)}(0)\big)
-\int_0^t ({\mathcal L}_K f)\big(Y^{(\infty,K)}(s)\big)\,\mathrm{d}s,\quad t\ge0
\end{equation}
is a martingale with respect to ${\mathcal F}^{(\infty,K)}(t)$, $t\ge0$ for any twice continuously differentiable function $f:\,\rr^K\to\rr$ which is bounded together with all its first and second order derivatives. To this end, we use the martingale property of the process
\begin{equation*}
f\big(Y^{(\infty,K)}(t)\big)-f\big(Y^{(\infty,K)}(0)\big)
-\int_0^t ({\mathcal L}_\infty f)\big(Y^{(\infty)}(s)\big)\,\mathrm{d}s,\quad t\ge0
\end{equation*}
with respect to ${\mathcal F}^{(\infty)}(t)$, $t\ge0$ (due to It\^o's formula) to find
\begin{eqnarray*}
\ev\big[f\big(Y^{(\infty,K)}(t)\big)\big|{\mathcal F}^{(\infty,K)}(s)\big]
&=&\ev\big[\ev\big[f\big(Y^{(\infty,K)}(t)\big)\big|{\mathcal F}^{(\infty)}(s)\big]\big|{\mathcal F}^{(\infty,K)}(s)\big] \\
&=& f\big(Y^{(\infty,K)}(s)\big)
+\ev\bigg[\ev\bigg[\int_s^t ({\mathcal L}_\infty f)\big(Y^{(\infty)}(u)\big)\,\mathrm{d}u\bigg|{\mathcal F}^{(\infty)}(s)\bigg]\bigg|{\mathcal F}^{(\infty,K)}(s)\bigg] \\
&=& f\big(Y^{(\infty,K)}(s)\big)+\ev\bigg[\int_s^t ({\mathcal L}_K f)\big(Y^{(\infty,K)}(u)\big)\,\mathrm{d}u
\bigg|{\mathcal F}^{(\infty,K)}(s)\bigg] \\
&&\quad\quad\quad\quad\quad\quad\;
+\ev\bigg[\int_s^t \big(({\mathcal L}_\infty-{\mathcal L}_K)f\big)\big(Y^{(\infty)}(u)\big)\,\mathrm{d}u
\bigg|{\mathcal F}^{(\infty,K)}(s)\bigg]
\end{eqnarray*}
for any $t\ge s\ge0$. It therefore remains to show that the last summand in the latter expression is equal to zero. In addition, thanks to Fubini's Theorem this statement can be further reduced to
\begin{equation}\label{prop_claim1'}
\ev\Big[\big(({\mathcal L}_\infty-{\mathcal L}_K)f\big)\big(Y^{(\infty)}(u)\big)\Big|{\mathcal F}^{(\infty,K)}(s)\Big]=0,\quad u\ge s\ge 0. 
\end{equation}

\smallskip

We first establish \eqref{prop_claim1'} assuming the following claim and then supply the proof of the claim in Step 2 below. 

\medskip

\noindent\textbf{Claim 2.} For any $u\ge 0$, the $\sigma$-algebra ${\mathcal F}^{(\infty,K)}(u)$ and the random vector $\big(X_K(u)-X_{K+1}(u),\,X_{K+1}(u)-X_{K+2}(u),\,\ldots\big)$ are independent. 

\medskip

Given the claim we compute
\begin{eqnarray*}
\ev\Big[\big(({\mathcal L}_\infty-{\mathcal L}_K)f\big)\big(Y^{(\infty)}(u)\big)\Big|{\mathcal F}^{(\infty,K)}(s)\Big] 
&=&\ev\Big[\ev\Big[\big(({\mathcal L}_\infty-{\mathcal L}_K)f\big)\big(Y^{(\infty)}(u)\big)
\Big|{\mathcal F}^{(\infty,K)}(u)\Big]\Big|{\mathcal F}^{(\infty,K)}(s)\Big] \\
&=&\ev\Big[\ev\Big[\big(({\mathcal L}_\infty-{\mathcal L}_K)f\big)\big(Y^{(\infty)}(u)\big)
\Big|Y^{(\infty,K)}(u)\Big]\Big|{\mathcal F}^{(\infty,K)}(s)\Big]. 
\end{eqnarray*}
Hence, to obtain \eqref{prop_claim1'} it suffices to prove
\begin{equation}\label{prop_claim2'}
\ev\Big[\big(({\mathcal L}_\infty-{\mathcal L}_K)f\big)\big(Y^{(\infty)}(u)\big)\Big|Y^{(\infty,K)}(u)\Big]=0,
\quad u\ge0.
\end{equation}
At this point, we use formula \eqref{1D_inf} to write
\begin{eqnarray*}
&&\ev\Big[\big(({\mathcal L}_\infty-{\mathcal L}_K)f\big)\big(Y^{(\infty)}(u)\big)\Big|Y^{(\infty,K)}(u)\Big] \\
&&=\sum_{l=K}^{K+d-1} \int_{\rr^d} 
\bigg(\frac{\partial f}{\partial x_1}\big(Y^{(\infty,K)}(u)\big)\,r_{l1}
+\sum_{k=1}^{K-1} \frac{\partial f}{\partial y_k}\big(Y^{(\infty,K)}(u)\big)
\,\big(r_{lk}-r_{l(k+1)}\big)\bigg) \\
&&\qquad\qquad\qquad\quad\quad\;\;
\cdot\big(U'(z_l)-\nu_l\big)\,\prod_{m=K}^{K+d-1} \frac{1}{{\mathcal Z}_m}\,\exp(2\,U(z_m)-2\,\nu_m\,z_m)\,\mathrm{d}z_m.
\end{eqnarray*} 
Note that each of the $d$ integrands is a first order partial derivative of a function vanishing at infinity, so that all $d$ integrals are equal to zero. The theorem now readily follows. 

\medskip

\noindent\textbf{Step 2.} It remains to supply the proof of Claim 2. To this end, we will analyze the time-reversed process $\widehat{Y}^{(\infty,K)}(s):=Y^{(\infty,K)}(u-s)$, $s\in[0,u]$ in the filtration $\widehat{\mathcal{F}}^{(\infty)}(s)$, $s\in[0,u]$ generated by $Y^{(\infty)}(u-s)$, $s\in[0,u]$ and verify that its dynamics does not depend on $\big(X_K(u)-X_{K+1}(u),\,X_{K+1}(u)-X_{K+2}(u),\,\ldots\big)$.

\medskip

Recalling the notations $A^{(K)}$ for the covariance matrix $(a_{kl})_{k,l=1}^K$ and $\|\cdot\|_{(A^{(K)})^{-1}}$ for the Euclidean norm associated with its inverse we find
\begin{equation}\label{twice_ent}
\ev\bigg[\int_0^u \Big\|\Big(\mu_k+\sum_{l=k}^{k+d-1} U'\big(X_l(s)-X_{l+1}(s)\big)\,r_{lk}\Big)_{1\le k\le K}\Big\|^2_{(A^{(K)})^{-1}}\,\mathrm{d}s\bigg]<\infty
\end{equation}
thanks to Fubini's Theorem, \eqref{1D_inf}, and \eqref{Fisher'}. Therefore, according to \cite[Theorem 7.2]{LS}, the law of the process $\big(X_1(s),\,X_2(s),\,\ldots,\,X_K(s)\big)$, $s\in[0,u]$ is absolutely continuous with respect to the law of a $K$-dimensional Brownian motion with zero drift vector, covariance matrix $A^{(K)}$ and initial condition $\big(X_1(0),\,X_2(0),\,\ldots,\,X_K(0)\big)$. Moreover, by \cite[Theorem 7.6, equation (7.30)]{LS} the relative entropy corresponding to this change of measure is given by one half of the quantity in \eqref{twice_ent} and is, hence, finite. In other words, the relative entropy of the law of the process $Y^{(\infty,K)}(s)$, $s\in[0,u]$ with respect to the law of a Brownian motion with zero drift vector, covariance matrix $\Theta^{(K)}$ (introduced in Step 4 of the proof of Proposition \ref{prop_cons}) and initial condition $Y^{(\infty,K)}(0)$ is finite. At this point, as in Step 4 of the proof of Proposition \ref{prop_cons}, one concludes that the law of the time-reversed process $\widehat{Y}^{(\infty,K)}(s)$, $s\in[0,u]$ has finite relative entropy with respect to the law of a Brownian motion on the time interval $[0,u]$ with zero drift vector and covariance matrix $\Theta^{(K)}$, conditioned to have the distribution of $Y^{(\infty,K)}(0)$ at time $u$ (that is, the process of \eqref{gamma2_eq} with $K$ replaced by $(K+1)$).

\medskip

Recall that in Step 4 of the proof of Proposition \ref{prop_cons} it has been established that the law of the process of \eqref{gamma2_eq} is locally absolutely continuous with respect to the law of a Brownian motion $\widehat{b}$ on the time interval $[0,u)$ with zero drift vector and covariance matrix $\Theta^{(K)}$. Consequently the law of the process $\widehat{Y}^{(\infty,K)}(s)$, $s\in[0,u)$ is also locally absolutely continuous with respect to the law of such a Brownian motion, so that $\widehat{Y}^{(\infty,K)}$ must satisfy an SDE of the form
\begin{equation}\label{gamma_eq'}
\mathrm{d}\widehat{Y}^{(\infty,K)}(s)=\gamma(s)\,\mathrm{d}s+\mathrm{d}\widehat{b}(s),\quad s\in[0,u)
\end{equation}
by \cite[Theorem 7.11]{LS}. Here $\gamma$ is a non-anticipative process in the filtration $\widehat{\mathcal F}^{(\infty)}(s)$, $s\in[0,u)$. Moreover, according to \cite[Theorems 7.5, 7.6]{LS}, the densities associated with the latter absolute continuity relation are given by
\begin{equation}\label{Gir2'}
\exp\bigg(\int_0^s \gamma(q)^T\big(\Theta^{(K)}\big)^{-1}\mathrm{d}\widehat{Y}^{(\infty,K)}(q)
-\frac{1}{2}\int_0^s \|\gamma(q)\|^2_{(\Theta^{(K)})^{-1}}\,\mathrm{d}q\bigg),\quad s\in[0,u).
\end{equation}

\smallskip

To determine $\gamma$ we consider the decomposition $\gamma(s)=:\gamma^{(1)}(s)+\gamma^{(2)}(s,\widehat{Y}^{(\infty,K)}(s))$, $s\in[0,u)$ where $\gamma^{(2)}$ is the same as in \eqref{gamma2_eq}, but with $(K-1)$ replaced by $K$. Combining \eqref{Gir2'} with the expression for the change of measure densities in \eqref{Gir1} (with $(K-1)$ replaced by $K$) and Girsanov's Theorem in the form of \cite[Theorem 3.5.1]{KS} we conclude that the densities associated with the change of measure from \eqref{gamma2_eq} (with $(K-1)$ replaced by $K$) to \eqref{gamma_eq'} are given by
\begin{equation*}
\exp\bigg(\int_0^s \gamma^{(1)}(q)\big(\Theta^{(K)}\big)^{-1}\mathrm{d}\widehat{b}(q)
+\frac{1}{2}\int_0^s \big\|\gamma^{(1)}(q)\big\|^2_{(\Theta^{(K)})^{-1}}\,\mathrm{d}q\bigg),\quad s\in[0,u).
\end{equation*}
Now, the finiteness of the corresponding relative entropy together with the simple localization argument in \cite[proof of Lemma 2.6]{Fo1} give
\begin{equation*}
\ev\bigg[\int_0^s \big\|\gamma^{(1)}(q)\big\|^2_{(\Theta^{(K)})^{-1}}\,\mathrm{d}q\bigg]<\infty,\quad s\in[0,u),
\end{equation*}
so that
\begin{equation*}
\ev\bigg[\int_0^s \big\|\gamma^{(1)}(q)\big\|^2\,\mathrm{d}q\bigg]<\infty,\quad s\in[0,u)
\end{equation*}
where $\|\cdot\|$ is the standard Euclidean norm on $\rr^K$. This estimate yields further
\begin{equation}\label{a_priori'}
\ev\bigg[\int_0^s \big\|\gamma(q)\big\|^2\,\mathrm{d}q\bigg]<\infty,\quad s\in[0,u)
\end{equation}
thanks to $\gamma(q)=\gamma^{(1)}(q)+\gamma^{(2)}(q,\widehat{Y}^{(\infty,K)}(q))$, the uniform sub-linear growth of $\gamma^{(2)}(q,\widehat{Y}^{(\infty,K)}(q))$ in $\widehat{Y}^{(\infty,K)}(q)$ on compact subsets of $q\in[0,u)$ (see the paragraph following \eqref{gamma2_eq}), \eqref{1D_inf}, and the finite second moment condition in Assumption \ref{MP_asmp}. 
 
\medskip

In view of the estimate \eqref{a_priori'} one can now apply \cite[Proposition 2.5]{Fo2} (more precisely its straightforward extension to the multidimensional setting at hand) to obtain the representation
\begin{equation*}
\gamma(s)=\lim_{\epsilon\downarrow0} \;
\epsilon^{-1}\;\ev\big[\widehat{Y}^{(\infty,K)}(s+\epsilon)-\widehat{Y}^{(\infty,K)}(s)
\big|\widehat{\mathcal F}^{(\infty)}(s)\big]
=\lim_{\epsilon\downarrow0} \;\epsilon^{-1}\;\int_{s-\epsilon}^s \gamma(q)\,\mathrm{d}q
\end{equation*}
for Lebesgue almost every $s\in[0,u)$ where the limits should be understood in the $L^2$ sense. In view of the assumed Markov property of $X$ (and, hence, of $Y^{(\infty)}$) the latter representation reduces to
\begin{equation}\label{gamma_char'}
\gamma(s)=\lim_{\epsilon\downarrow0} \;
\epsilon^{-1}\;\ev^\pp\big[\widehat{Y}^{(\infty,K)}(s+\epsilon)-\widehat{Y}^{(\infty,K)}(s)\big|Y^{(\infty)}(u-s)\big]
=\lim_{\epsilon\downarrow0} \;\epsilon^{-1}\;\int_{s-\epsilon}^s \gamma(q)\,\mathrm{d}q.
\end{equation}

\smallskip

Lastly, pick an $L\ge K+d$ and a function $f:\,\rr^\nn\to\rr$ depending only in the first $L$ coordinates which is twice continuously differentiable and has compact support with respect to these coordinates. Then, on the one hand,
\begin{equation}\label{gamma1'}
\lim_{\epsilon\downarrow0}\;\epsilon^{-1}\;\ev\Big[\big(\widehat{Y}^{(\infty,K)}(s+\epsilon)
-\widehat{Y}^{(\infty,K)}(s)\big) f\big(Y^{(\infty)}(u-s)\big)\Big]=\ev\big[\gamma(s)\,f\big(Y^{(\infty)}(u-s)\big)\big]
\end{equation}
for Lebesgue almost every $s\in[0,u)$ by virtue of \eqref{gamma_char'}. On the other hand, with the notations $D^{(\infty,K)}$ for the drift coefficient of $Y^{(\infty,K)}$ with respect to the filtration ${\mathcal F}^{(\infty)}(t)$, $t\ge0$ and 
\begin{eqnarray*}
{\mathcal L}_{\infty,L}&:=&
\frac{1}{2}\,a_{11}\,\frac{\partial^2}{\partial x_1^2}
+\sum_{k=1}^{L-1} (a_{1k}-a_{1(k+1)})\,\frac{\partial^2}{\partial x_1\partial y_k}
+\mu_1\,\frac{\partial}{\partial x_1}
+\sum_{l=1}^\infty U'(y_l)\,r_{l1}\,\frac{\partial}{\partial x_1} \\
&&+\frac{1}{2}\,\sum_{k,l=1}^{L-1} \big(a_{kl}+a_{(k+1)(l+1)}-a_{(k+1)l}-a_{k(l+1)}\big)
\,\frac{\partial^2}{\partial y_k\,\partial y_l}
+\sum_{k=1}^{L-1} \big(\mu_k-\mu_{k+1}\big)\,\frac{\partial}{\partial y_k} \\
&&+\sum_{k=1}^{L-1} \sum_{l=k}^\infty U'(y_l)\,\big(r_{lk}-r_{l(k+1)}\big)\,\frac{\partial}{\partial y_k}\,, 
\end{eqnarray*}
one computes
\begin{equation*}
\begin{split}
&\lim_{\epsilon\downarrow0}\;\epsilon^{-1}\;\ev\Big[\big(\widehat{Y}^{(\infty,K)}(s+\epsilon)
-\widehat{Y}^{(\infty,K)}(s)\big) f\big(Y^{(\infty)}(u-s)\big)\Big] \\
&=-\lim_{\epsilon\downarrow0}\;\epsilon^{-1}\;\ev\Big[Y^{(\infty,K)}(u-s)\,f\big(Y^{(\infty)}(u-s)\big)
-Y^{(\infty,K)}(u-s-\epsilon)\,f\big(Y^{(\infty)}(u-s-\epsilon)\big)\Big] \\
&\quad+\lim_{\epsilon\downarrow0}\;\epsilon^{-1}\;\ev\Big[Y^{(\infty,K)}(u-s-\epsilon)
\Big(f\big(Y^{(\infty)}(u-s)\big)-f\big(Y^{(\infty)}(u-s-\epsilon)\big)\Big)\Big] \\
&=-\ev\Big[Y^{(\infty,K)}(u-s)\,\big({\mathcal L}_{\infty,L} f\big)\big(Y^{(\infty)}(u-s)\big)\Big]
-\ev\Big[f\big(Y^{(\infty)}(u-s)\,D^{(\infty,K)}(u-s)\Big] \\
&\quad-\ev\bigg[\bigg(\frac{\partial f}{\partial x_1}\big(Y^{(\infty)}(u-s)\big)\,\Theta^{(L)}_{1k}
+\sum_{l=1}^{L-1} \frac{\partial f}{\partial y_l}\big(Y^{(\infty)}(u-s)\big)\,\Theta^{(L)}_{(l+1)k}\bigg)_{1\le k\le K}\bigg] \\
&\quad+\ev\Big[Y^{(\infty,K)}(u-s)\,\big({\mathcal L}_{\infty,L} f\big)\big(Y^{(\infty)}(u-s)\big)\Big] 
\end{split}
\end{equation*}
for Lebesgue almost every $s\in[0,u)$. Here we have used It\^o's formula and forward in time versions of \eqref{gamma_char'}. By cancelling out the two identical expectations and integrating by parts one can simplify the latter expression to
\begin{equation}\label{gamma2'}
\begin{split}
\ev\bigg[f\big(Y^{(\infty)}(u-s)\big)\bigg(-D^{(\infty,K)}(u-s)
+\bigg(-\frac{Y^{(\infty)}_1(u-s)-\big(\mu_1+\sum_{l=1}^d r_{l1}\nu_l\big)(u-s)}{a_{11}(u-s)}\,\Theta^{(L)}_{1k} \qquad\qquad\quad\\
+\sum_{l=2}^L \Big(2\,U'\big(Y^{(\infty)}_l(u-s)\big)-2\,\nu_{l-1}\Big)\,\Theta^{(L)}_{lk}\bigg)_{1\le k\le K}\bigg)\bigg].
\end{split}
\end{equation}
Since $\gamma(s)$ is measurable with respect to the $\sigma$-algebra generated by $Y^{(\infty)}(u-s)$, a comparison of \eqref{gamma1'} and \eqref{gamma2'} allows to determine $\gamma(s)$ for Lebesgue almost every $s\in[0,u)$. Using Assumption \ref{OO_asmp} one can deduce further that, for such $s$, the drift coefficient $\gamma(s)$ is a deterministic function of $Y^{(\infty,K)}(u-s)=\widehat{Y}^{(\infty,K)}(s)$ only. The proof of weak uniqueness in part (a) of Proposition \ref{prop_cons} can be now repeated mutatis mutandis to show weak uniqueness for the SDE resulting from \eqref{gamma_eq'}, and one concludes that the dynamics of the process $\widehat{Y}^{(\infty,K)}(s)$, $s\in[0,u]$ in the filtration $\widehat{\mathcal F}^{(\infty)}(s)$, $s\in[0,u]$ does not depend on $\big(X_K(u)-X_{K+1}(u),\,X_{K+1}(u)-X_{K+2}(u),\,\ldots\big)$. This gives Claim 2. \ep

\section{Examples}\label{sec_examples}

To illustrate the wide applicability of our results we give two examples of infinite Brownian particle systems arising in the area of integrable probability to which they apply. 

\begin{example}[$\beta$ analogues of Brownian TASEP] 
Fix a $\beta>0$ and consider the $\beta$ corners process which is the probability measure on the Gelfand-Tseitlin cone
\begin{equation*}
\big\{x=(x^j_i:\,1\le i\le j\le N)\in\rr^{N(N+1)/2}:\;x^{j+1}_i\le x^j_i\le x^{j+1}_{i+1},\,1\le i\le j\le N-1\big\}
\end{equation*} 
with density proportional to
\begin{equation*}
\prod_{1\le i<j\le N} (x^N_j-x^N_i)\prod_{i=1}^N \exp\Big(-\frac{(x_i^N)^2}{2}\Big)\prod_{j=1}^N \prod_{1\le i<i'\le j} 
(x^j_i-x^j_{i'})^{2-\beta}\prod_{\iota=1}^j \prod_{\iota'=1}^{j+1} |x^j_\iota-x^{j+1}_{\iota'}|^{\beta/2-1}.
\end{equation*}
The interest in this measure comes from the fact that, for $\beta=1$, $2$, and $4$, it describes the joint distribution of the eigenvalues of the $1\times1,\,2\times2,\,\ldots,\,N\times N$ top left corners of a random matrix from the Gaussian orthogonal ensemble (GOE), Gaussian unitary ensemble (GUE), and Gaussian symplectic ensemble (GSE), respectively (see e.g. \cite{Ne} and \cite{GS1} for more details). In the papers \cite{Wa}, \cite{GS1} dynamic versions of the $\beta$ corners process were introduced for $\beta=2$ and $\beta>2$, respectively, with the latter being given by the diffusion process
\begin{equation*}
\mathrm{d}R^j_i(t)=\mathrm{d}w^j_i(t)-\sum_{i'\neq i} \frac{\beta/2-1}{R^j_i(t)-R^j_{i'}(t)}\,\mathrm{d}t
+\sum_{i'=1}^{j-1} \frac{\beta/2-1}{R^j_i(t)-R^{j-1}_{i'}(t)}\,\mathrm{d}t,\quad 1\le i\le j\le N
\end{equation*}
on the Gelftand-Tseitlin cone. Here $w^j_i$, $1\le i\le j\le N$ are independent standard Brownian motions. 

\medskip

The study of the extremal particles $R^N_N,\,R^{N-1}_{N-1},\,\ldots$ is of particular interest due to their relations with objects in the Kardar-Parisi-Zhang (KPZ) universality class. Indeed, for $\beta=2$, the extremal particles form a Brownian analogue of the totally asymmetric simple exclusion process (see \cite{Wa} for more details) which is known to belong to the KPZ universality class (see e.g. \cite{BF} and the references therein). For $\beta>2$ the situation is more intricate, but one can show that, at least for $\beta\ge4$, the process of extremal particles $R^N_N,\,R^{N-1}_{N-1},\,\ldots$, seen from the point of view of the particle $R^N_N$, converges (after the appropriate rescaling) in the limit $N\to\infty$ to the spacings in the infinite Brownian particle system
\begin{equation}\label{beta_TASEP}
\mathrm{d}X_k(t)=\frac{1}{\sqrt{2}}\,\mathrm{d}B_k(t)+\frac{\beta/4-1/2}{X_k(t)-X_{k+1}(t)}\,\mathrm{d}t,\quad k\in\nn
\end{equation}
endowed with the quasi-stationary initial condition 
\begin{equation}\label{beta_TASEP_ic}
\delta_0(\mathrm{d}x_1)\,\prod_{k=1}^\infty \frac{1}{\Gamma(\beta/2)}\,\Big(\frac{\beta}{2}\Big)^{\beta/2}
\,(x_k-x_{k+1})^{\beta/2-1}\,\exp\Big(-\frac{\beta}{2}\,(x_k-x_{k+1})\Big)\,\mathrm{d}(x_k-x_{k+1}).
\end{equation}
Hereby, $\Gamma(\cdot)$ is the Gamma function. We refer to \cite{GS2} for more details and note that in \eqref{beta_TASEP} time is slowed down by a factor of two compared to the setting there to concur with the normalizations of Assumption \ref{OO_asmp}. 

\medskip

One can put the system \eqref{beta_TASEP} into the framework of \eqref{main_sde} by setting $\mu_k=\mu$, $k\in\nn$ for an arbitrary $\mu>0$, $U(z)=\big(\frac{\beta}{4}-\frac{1}{2}\big)\log z-\mu z$, and $r_{lk}=\mathbf{1}_{k=l}$, $(k,l)\in\nn^2$. From Remark \ref{OO_rmk} it follows directly that Assumption \ref{OO_asmp} is satisfied. Moreover, the density $z^{\beta/2-1}\,e^{-\mu z}\,\mathbf{1}_{(0,\infty)}(z)$ can be normalized to the probability density of the appropriate Gamma distribution (with $\mu=\frac{\beta}{2}$ corresponding to \eqref{beta_TASEP_ic}), the latter has a finite second moment, and its Fisher information is finite iff
\begin{equation*}
\int_0^\infty z^{\beta/2-1}e^{-\mu z}\,\frac{1}{z^2}\,\mathrm{d}z<\infty\,,
\end{equation*}   
that is, iff $\beta>4$. It follows that, for $\beta>4$, the results of Theorems \ref{main_thm1} and \ref{main_thm2} apply to the system \eqref{beta_TASEP} and, in particular, identify \textit{all} measures 
\begin{equation}
\delta_0(\mathrm{d}x_1)\,\prod_{k=1}^\infty \frac{1}{\Gamma(\beta/2)}\,\mu^{\beta/2}
\,(x_k-x_{k+1})^{\beta/2-1}\,\exp\Big(-\mu(x_k-x_{k+1})\Big)\,\mathrm{d}(x_k-x_{k+1}),\quad\mu>0
\end{equation}
as quasi-stationary for that system.  
\end{example}

\begin{example}[O'Connell-Yor semi-discrete polymer \& Brownian $q$-TASEP] 
Fix an $N\in\nn$ and a $t\ge0$, and imagine a right and up polymer path in the plane connecting $(0,0)$ to $(t,N)$ by following the line $\rr\times\{0\}$ up to the point $(t_1,0)$ for some $0\le t_1\le t$, then moving up to $(t_1,1)$ and following the line $\rr\times\{1\}$ up to the point $(t_2,1)$ for some $t_1\le t_2\le t$, then moving up to $(t_2,2)$ etc. until the point $(t,N)$ is reached. In the O'Connell-Yor semi-discrete polymer model (see \cite{OY}, \cite{MO}, \cite{OC}) every such path is assigned a weight proportional to
\begin{equation*}
\exp\big(b_0(t_1)+(b_1(t_2)-b_1(t_1))+\ldots+(b_N(t)-b_N(t_N))\big)
\end{equation*}
where $b_0,\,b_1,\,\ldots,\,b_N$ are independent standard Brownian motions. The logarithmic partition function associated with this measure is then given by
\begin{equation*}
Z_N(t):=\log\int_{0\le t_1\le t_2\le\cdots\le t_N\le t} 
\exp\big(b_0(t_1)+(b_1(t_2)-b_1(t_1))+\ldots+(b_N(t)-b_N(t_N))\big)
\,\mathrm{d}t_1\,\mathrm{d}t_2\,\ldots\,\mathrm{d}t_N.
\end{equation*}

\smallskip

As one lets $N\in\nn$ and $t\ge0$ vary, one finds (see \cite{OC}) that each process $\big(Z_0(t),\,Z_1(t),\,\ldots,\,Z_N(t)\big)$, $t\ge0$ is a diffusion satisfying
\begin{equation}\label{OCY}
\begin{split}
&\mathrm{d}Z_k(t)=\mathrm{d}b_k(t)+e^{Z_{k-1}(t)-Z_k(t)}\,\mathrm{d}t,\quad k=1,\,2,\,\ldots,\,N,\\ 
&\mathrm{d}Z_0(t)=\mathrm{d}b_0(t).
\end{split}
\end{equation}
The same diffusion also appears in the context of the Whittaker growth model (see \cite[Definition 4.1.26]{BC}) and provides the scaling limit of the $q$-TASEP process introducted in \cite[Section 3.3.2]{BC} (see \cite[Remark 4.1.28]{BC} for more details). It can be also viewed as a series of queues in tandem as explained in \cite{OY}. A natural infinite-dimensional analogue of \eqref{OCY} is the system of SDEs
\begin{equation}\label{OCY1}
\mathrm{d}X_k(t)=\frac{1}{\sqrt{2}}\,\mathrm{d}b_k(t)+\frac{1}{2}\,e^{-(X_k(t)-X_{k+1}(t))}\,\mathrm{d}t,\quad k=1,\,2,\,\ldots
\end{equation}
which, in particular, has an interpretation as an infinite series of queues in tandem similar to the one for the system \eqref{OCY}. As in the previous example the slowdown by a factor of two is solely for consistency with Assumption \ref{OO_asmp}. 

\medskip

The system \eqref{OCY1} falls into the framework of \eqref{main_sde} with $\mu_k=\frac{\mu}{2}$, $k\in\nn$ for an arbitrary $\mu>0$, $U(z)=-\frac{1}{2}(\mu z+e^{-z})$, and $r_{lk}=\mathbf{1}_{k=l}$, $(k,l)\in\nn^2$. According to Remark \ref{OO_rmk} the system \eqref{OCY1} satisfies Assumption \ref{OO_asmp}. In addition, the density $\exp\big(-\mu z-e^{-z}\big)$ can be normalized to a probability measure which has finite second moment and Fisher information, so that Assumption \ref{MP_asmp} is fulfilled. Thus, Theorems \ref{main_thm1} and \ref{main_thm2} apply to the system \eqref{OCY1} and, in particular, show that \textit{all} measures of the form 
\begin{equation}
\delta_0(\mathrm{d}x_1)\,\prod_{k=1}^\infty \frac{1}{\Gamma(\mu)}\,\exp\Big(-\mu(x_k-x_{k+1})-e^{-(x_k-x_{k+1})}\Big)\,\mathrm{d}(x_k-x_{k+1}),\quad\mu>0
\end{equation}
are quasi-stationary for that system. As before $\Gamma(\cdot)$ stands for the Gamma function.
\end{example}

\bigskip\bigskip

\end{document}